\documentclass{amsart}
\usepackage{graphicx}
\newtheorem{theorem}{Theorem}[section]
\newtheorem{corollary}[theorem]{Corollary}
\newtheorem{lemma}[theorem]{Lemma}

\theoremstyle{definition}
\newtheorem{definition}[theorem]{Definition}
\theoremstyle{remark}

\numberwithin{equation}{section}

\newcommand{\BB}{\mathbb B}
\newcommand{\CC}{\mathbb C}
\newcommand{\RR}{\mathbb R}

\def\cC{{\mathcal C}}

\begin{document}

\title{A Generalization of Forelli's Theorem}%
\author{Jae-Cheon Joo, Kang-Tae Kim and Gerd Schmalz}%
\address{(Joo) Department of Mathematics, Pusan National University, Busan 609-735 The Republic of Korea}%
\email{jcjoo91@pusan.ac.kr}
\address{(Kim) Department of Mathematics,
Pohang University of Science and Technology, Pohang
790-784 The Republic of Korea}%
\email{kimkt@postech.ac.kr}%
\address{(Schmalz) School of Science and Technology, University of New England, Armidale, NSW 2351, Australia}%
\email{gerd@turing.une.edu.au}%


\begin{abstract}
The purpose of this paper is to present a generalization of Forelli's theorem. In particular, we prove an all dimensional version of the two-dimensional theorem of Chirka \cite{Chirka} of 2005.
\end{abstract}
\maketitle

\section{Introduction}

Classical Forelli's theorem (\cite{Forelli}; see also \cite{Shabat}, \cite{Stoll}) states

\begin{theorem}[Forelli]\label{forelli}
Let $f\colon \BB^n \to \CC$ be a function. If it satisfies the
following two conditions:
\begin{itemize}
\item[(\romannumeral 1)]  $f \in \cC^\infty (0)$,
\end{itemize}
and %
\begin{itemize}
\item[(\romannumeral 2)] for every unit vector $v=(v_1, \ldots,
v_n) \in \CC$, $f (\zeta v_1, \ldots, \zeta v_n)$ is holomorphic
in the single complex variable $\zeta$ with $|\zeta|<1$,
\end{itemize}
then $f$ is holomorphic.
\end{theorem}

The definition of the notation $f \in \cC^\infty(0)$ in the statement is as follows: for any positive integer $k$, there exists a polynomial  $p_k$ such that $f(z)-p_k(z)= \operatorname{o}(|z|^k)$.
\medskip

After a long period of almost no results, the following two generalizations have been presented:

\begin{theorem}[Chirka \cite{Chirka}]\label{chirka}
Let $\{S_\tau\}$ be a foliation of a domain $\Omega$ in a punctured ball $\mathbb B \setminus 0$ in ${\mathbb C}^2$ by holomorphic curves that are closed and smooth in $\mathbb B$, pass through the origin, and are pairwise transversal at $0$. Let $f$ be a function at $\mathbb B$ such that $f\in C^\infty(0)$ and all restrictions $f|_{S_\tau}$ are holomorphic. Then $f$ is holomorphic in $\Omega$ and, if $\Omega= \mathbb B\setminus 0$, then $f$ is holomorphic in $\mathbb B$.
\end{theorem}

\begin{theorem}[Kim-Poletsky-Schmalz \cite{Kim-Poletsky-Schmalz}]\label{KPS}
If $n\ge 1$ is an integer, and if $f\colon \mathbb B^n \to \mathbb C$ is a function with $f\in \mathcal C^\infty(0)$, which is annihilated by
$$
\bar{X} = \sum_{j=1}^n \alpha_j \bar{z}_j \frac{\partial}{\partial \bar{z}_j},
$$
where $\alpha_1,\dots,\alpha_n$ are real numbers, then $f$ is holomorphic on $\mathbb B^n$.
\end{theorem}

In both versions holomorphicity along straight lines through $0$ has been replaced by holomorphicity along some more general family of complex curves: in case of Theorem \ref{chirka} they are assumed to intersect transversely. On the other hand, in the case of Theorem \ref{KPS} the family of complex curves is generated by a holomorphic vector field and its integral curves. Thus, the two theorems complement each other in a sense; the general foliation considered by Chirka may not in general be generated by a contracting holomorphic vector field, whereas the leaves of foliation considered by Kim-Poletsky-Schmalz do not have to intersect mutually transversely at the origin (even after the re-parametrization of the leaves so that they intersect at the origin).
\medskip

The purpose of this paper is to present Theorem 2.1, an all-dimensional version of Theorem \ref{chirka}. This in particular answers the question posed by Chirka in \cite{Chirka}, p.\ 219.

\section{A generalization of Forelli's theorem in $\CC^n$}

First we define the concept of the local $\cC^k$ singular foliation by holomorphic curves. Let $\Delta\subset \mathbb C$ the unit disc and $S^{2n-1}$ the unit sphere in $\CC^n$ defined by the equation $|z_1|^2+\cdots+|z_n|^2=1$.

\begin{definition} \rm
Let $\ell$ be a positive integer. For a point $p \in \CC^n$, a \em local $\cC^\ell$ singular foliation at $p$ by holomorphic discs \rm is a $\cC^\ell$ map $h\colon \Delta \times S^{2n-1} \to \CC^n$ satisfying the following properties:
\begin{itemize}
\item[(1)] For each $v \in S^{2n-1}$, the correspondence $h(\cdot, v) \colon z \in \Delta \to h(z, v) \in \CC^n$ is a holomorphic embedding.
\item[(2)] $h(0,v) = p$ for every $v \in S^{2n-1}$.
\item[(3)]The image $h(\Delta \times S^{2n-1})$ contains an open neighborhood of $p$ in $\CC^n$.
\item[(4)] For each $v \in S^{2n-1}$, there exists $r_v>0$ such that $\displaystyle{\frac{\partial h}{\partial z}\Big|_{(0,v)} = r_v v}$.
\item[(5)] $h(z, e^{i\theta} v) = h(e^{i\theta}z, v)$ for any $\theta \in \RR$ and $z \in \Delta$.
\end{itemize}
\end{definition}

Throughout this paper, we shall consider the case $\ell=1$ only.
\medskip

We shall consider, from here on, only the case when $p$ is the origin. This singular foliation provides a parametrization of an open neighborhood of the origin in $\CC^n$ by $\Delta \times \mathbb C\mathbb P^{n-1}$. One can always choose coordinates around the origin in $\mathbb C^n$ such that any given direction $v^0 \in S^{2n-1}$ becomes $[1:0\dots:0]\in\mathbb C\mathbb P^{n-1}$ and hence a neighborhood of $v^0$ in $\mathbb C\mathbb P^{n-1}$ can be identified as a neighborhood of $0$ in $\mathbb C^{n-1}$ with coordinates $\displaystyle{(\lambda^1, \ldots, \lambda^0_{n-1}) := \Big( \frac{v^0_2}{v^0_1}, \ldots, \frac{v^0_n}{v^0_1}\Big)}$.  Then $h(z,v^0)$ can be understood in this coordinate system as $h(z,0)$, and in a neighborhood of the corresponding holomorphic curve, we may assume without loss of generality that it has the equation $h(z,0)=(z,0)$. Then for $z$ and $\lambda$, with both $|z|$ and $\|\lambda\|$ sufficiently small, the curve $z \to h(z,\lambda)$ is represented by the expression
$$
\left\{
\begin{array}{ccl}
z & = & z \\
w_1 & = & k_1(z,\lambda_1, \ldots, \lambda_{n-1}) \\
& \vdots & \\
w_{n-1} & = & k_{n-1} (z,\lambda_1, \ldots, \lambda_{n-1}).
\end{array}
\right.
$$
We use the short-hand notation $w=k(z,\lambda)$. Then $k$ satisfies
\begin{enumerate}
\item $k(z,\lambda)$ is $\cC^1$ in $z,\lambda$ and holomorphic in $z$.
\item $k(0,\lambda)=0$ for any $\lambda$
\item $\displaystyle{\frac{\partial k(z,\lambda)}{\partial z}\Big|_{z=0}}=\lambda$ and hence $k(z,\lambda)= z \lambda + \operatorname{o}(|z|)$.
\end{enumerate}

Notice that we are using the standard complex coordinate system $(z, w_1,\ldots, w_{n-1})$ for $\CC^n$ here.

\begin{theorem} \label{Chirka}
If $h$ is a local singular foliation of a domain $\Omega$ in $\CC^n$ and $f\colon\Omega \to \CC$ is a function satisfying:
\begin{itemize}
\item[(A)] $f \in \cC^\infty(0)\cap \cC^1(\Omega)$; and
\item[(B)] for every leaf $h(\cdot, \lambda)$ the composition $f\circ h(\cdot,\lambda)$ is a holomorphic function,
\end{itemize}
then $f$ is a holomorphic function on the intersection of $\Omega$ and some neighborhood of the origin.
\end{theorem}

Notice that the statement of this theorem in the case of $n=2$ is weaker than what was presented in \cite{Chirka}.  On the other hand, our proof here is not only valid for all dimensions, but also, even in dimension 2, somewhat more straightforward.

\begin{proof}
Let us use the notation $\partial f/\partial \bar w := (\partial f/\partial \bar w_1, \ldots, \partial f/\partial \bar w_{n-1})$. Of course the goal here is to establish that $\partial f/\partial \bar w=0$ at $w=0$.

Let $ F (z,\lambda_1,\dots,\lambda_{n-1}) = f(h(z,\lambda_1,\dots,\lambda_{n-1})) = f(z, k(z,\lambda))$.
First we prove that $\partial F / \partial \lambda_j$, $\partial F / \partial\bar\lambda_j$ and $\partial k_m / \partial \lambda_j$ and $\partial k_m / \partial\bar\lambda_j$ at $\lambda =0$ are also holomorphic in $z$. In fact, let $G\colon \Delta\times U\to \mathbb C$ be a $\cC^1$ function  that is holomorphic with respect to $z\in \Delta$, like $F,k$. Then
$$
\int G(z,\lambda)\,dz\wedge d\phi = 0
$$
for any function in $\phi(z)\in\cC^\infty_0(\Delta)$. By differentiating with respect to the parameter $\lambda_j$ or $\bar{\lambda}_j$ we get
$$
\int_{\Delta} \frac{\partial}{\partial\lambda_j} G(z,\lambda)\,dz\wedge d\phi = 0, \quad \int_{\Delta} \frac{\partial}{\partial\bar{\lambda}_j} G(z,\lambda)\,dz\wedge d\phi = 0,
$$
which shows that $\frac{\partial G}{\partial\lambda_j}$ and $\frac{\partial G}{\partial\bar{\lambda}_j}$ are holomorphic in $z$.

Since $F$ is holomorphic in $z$ for every $\lambda$ and since $F$ is $\cC^1$-smooth, $\partial F / \partial \lambda$ and $\partial F / \partial\bar\lambda$ at $\lambda =0$ are also holomorphic in $z$. By the chain rule,
\begin{align*}
\left.\frac{\partial F}{\partial \lambda}\right|_{\lambda =0} &= A \left.\frac{\partial f}{\partial w}\right|_{w =0} + \bar{B}
\left.\frac{\partial f}{\partial \bar{w}}\right|_{w =0}\\
\left.\frac{\partial F}{\partial \bar{\lambda}}\right|_{\lambda =0} &= B \left.\frac{\partial f}{\partial w}\right|_{w =0} + \bar{A} \left.\frac{\partial f}{\partial \bar{w}}\right|_{w =0},
\end{align*}
where  $\frac{\partial F}{\partial \lambda}$ denotes the column of partial derivatives $(\frac{\partial F}{\partial \lambda_1},\dots, \frac{\partial F}{\partial \lambda_{n-1}})$ etc.\ and $A=(A_{mj})$ and $B=(B_{mj})$ are the matrices whose entries are defined by the partial derivatives as follows:
 $A_{mj} = \frac{\partial k_m}{\partial \lambda_j}\big|_{\lambda=0}$ and $B_{mj}= \frac{\partial k_m}{\partial \bar{\lambda}_j}\big|_{\lambda=0}$.

From $k_m(z,\lambda)= z\lambda_m + \operatorname{o}(|z|^2)$ we get
\begin{equation}\label{AB}
A_{mj}=\frac{\partial k_m}{\partial \lambda_j}= \delta_{mj} z + \operatorname{o}(|z|),\quad B_{mj} =\frac{\partial k_m}{\partial \bar{\lambda}_j}=  \operatorname{o}(|z|).
\end{equation}

Hence the matrix $\tilde{A}:=\frac{1}{z}A=\operatorname{id}+\operatorname{o}(1)$ is invertible in some neighbourhood of $0$. Denote $(\tilde{A})^{-1}=C$. Then the entries of $C$ are holomorphic functions. Let
$$
H(z):= BC \left.\frac{\partial F}{\partial \lambda}\right|_{\lambda =0} - z\left.\frac{\partial F}{\partial \bar\lambda}\right|_{\lambda =0}= (BC\bar{B} - z\bar{A}) \left.\frac{\partial f}{\partial \bar w}\right|_{w=0}.
$$
Then  $H(z)$ is holomorphic on $\Delta$.

Now we shall prove:  $H(z) = 0, \forall z \in \Delta$. In order to show this, we need the following lemma. Below, the notation $\mathbb C[[z_1,\dots,z_n, \bar{z}_1,\dots,\bar{z}_k]]$ represents the local ring of formal power series in the variables $z_1,\dots,z_n, \bar{z}_1,\dots,\bar{z}_k$ at the origin. Of course, the unique maximal ideal is the set of all formal power series without the constant term.

\begin{lemma}
Let $\alpha_1, \ldots, \alpha_m\in \mathbb C[[z_1,\dots,z_n, \bar{z}_1,\dots,\bar{z}_k]]$,  $\psi\in \mathbb C[[z_1,\dots,z_n]]$  and $\varphi_1, \ldots, \varphi_m\in \mathfrak M$, where $\mathfrak M$ is the maximal ideal of the local ring $\mathbb C[[\bar{z}_1,\dots,\bar{z}_k]]$. Then
$$
\psi = \alpha_1\varphi_1  + \ldots + \alpha_m\varphi_m
$$
implies $\psi=0$.
\end{lemma}

{\it Proof.} Assume $\psi\neq 0$ and let $\tilde{\psi}$ be the lowest degree non-vanishing polynomial in $\psi$. Then
$$\tilde{\psi} = \tilde{\alpha}_1\tilde{\varphi}_1  + \ldots + \tilde{\alpha}_m\tilde{\varphi}_m,$$
where $\tilde{\alpha}_1, \ldots, \tilde{\alpha}_m$ are certain polynomials in  $z_1,\dots,z_n, \bar{z}_1,\dots,\bar{z}_k$ and $\tilde{\varphi}_1, \ldots, \tilde{\varphi}_m$ are polynomials in $\bar{z}_1,\dots,\bar{z}_k$ of positive degree. Now, $\tilde{\psi}$ does not contain variables  $\bar{z}_1,\dots,\bar{z}_k$, whereas each monomial in
$\tilde{\alpha}_1\tilde{\varphi}_1  + \ldots + \tilde{\alpha}_m\tilde{\varphi}_m$ does contain such variables. This contradiction shows that $\psi=0$.\hfill $\Box$
\medskip

The following statement is then immediate:

\begin{corollary} \label{rudi-2}
Assume that $\varphi_1, \ldots, \varphi_m$ and $\alpha_1, \ldots, \alpha_m$ are complex-valued functions
defined on a domain $\Omega$ in the complex plane $\CC$, which
enjoy the properties:
\begin{itemize}
\item[(a)] $\varphi_k $ has a formal Taylor expansion at $p$, and
\item[(b)] $\alpha_k$ is conjugate-holomorphic in an open
neighborhood of $p$ with $\alpha_k(p)=0$
\end{itemize}
for $k=1,\ldots,m$.
If $\psi := \alpha_1\varphi_1  + \ldots + \alpha_m\varphi_m$ is holomorphic in $\Omega$, then $\psi$ is identically zero.
\end{corollary}
\medskip

Now we return to $H(z)$.
From (\ref{AB}) it follows that the anti-holomorphic terms
$\dfrac{\partial \bar{k}_m}{\partial \lambda_j}$ and $\dfrac{\partial
\bar{k}_m}{\partial \bar{\lambda}_j}$
vanish at $\zeta=0$.  Therefore, the components of $H(z)$ have the form of a function $\psi$ from the Corollary \ref{rudi-2} with the
$\varphi$'s being products of $\dfrac{\partial f}{\partial \bar{\lambda}_j}$ and some holomorphic factors from the matrices $A,B,C$  and the $\alpha$'s being products of some antiholomorphic factors from the matrices $\bar{A}, \bar{B}$. It follows $H(z)\equiv 0$ and hence
$$
(BC\bar{B} - z\bar{A}) \left.\frac{\partial f}{\partial \bar w}\right|_{w=0}\equiv 0.
$$
Finally, even though $BC\bar{B} - z\bar{A}$ vanishes at the origin, its determinant equals
$$
\det(BC\bar{B} - z\bar{A})=(-1)^{n-1} |z|^{2n-2} + \operatorname{o}(|z|^{2n-2})
$$
and therefore has no zeroes in some punctured neighborhood of $0$. It follows that
$$
\left.\frac{\partial f}{\partial \bar{w}}\right|_{w=0}=0.
$$
Since $\{w=0\}$ was an arbitrary line leaf the Cauchy-Riemann equations are satisfied transversally to each leaf.
This completes the proof.
\end{proof}

Finally, for the global version of generalized Forelli's theorem, we recall that the following definition of global singular foliation.

\begin{definition} \rm
Let $\Omega$ be a domain in $\CC^n$ containing the origin. By a \em $\cC^1$ singular foliation at $0$ by holomorphic discs \rm we mean a $\cC^1$ map $h\colon \Delta \times S^{2n-1} \to \Omega$ satisfying the following properties:
\begin{itemize}
\item[(1)] For each $v \in S^{2n-1}$, the correspondence $h(\cdot, v) \colon z \in \Delta \to h(z, v) \in \CC^n$ is a holomorphic embedding.
\item[(2)] $h(0,v) = 0$ for every $v \in S^{2n-1}$.
\item[(3)] $h(\Delta \times S^{2n-1})= \Omega$.
\item[(4)] For each $v \in S^{2n-1}$, there exists $r_v>0$ such that $\displaystyle{\frac{\partial h}{\partial z}\Big|_{(0,v)} = r_v v}$.
\item[(5)] $h(z, e^{i\theta} v) = h(e^{i\theta}z, v)$ for any $\theta \in \RR$ and $z \in \Delta$.
\end{itemize}
\end{definition}

Then we present

\begin{corollary}
Let $S_{\lambda}$ be the typical leaf (i.e., a holomorphic disc) of the singular foliation as above of a domain $\Omega\subset\mathbb C^n$ and $f$ a function that is $\cC^\infty(0) $ and $\cC^1(\Omega)$ such that the restrictions $f|{S_\lambda}$ are holomorphic. Then $f$ is holomorphic on $\Omega$
\end{corollary}

This follows from Theorem \ref{Chirka} and Chirka's curvilinear Hartogs' lemma from \cite{Chirka} (cf.\ \cite{Hartogs}).

\bigskip
\bf Acknowledgements. \rm Research of the second named author is supported in part by the Grant 4.0006570 (2011) of The Basic Science Research Institute of Pohang University of Science and Technology, The Republic of Korea.


\end{document}